

\magnification=\magstep1
\baselineskip=15pt
\parskip=4pt

\def\({\left(}
\def\){\right)}
\def\z{\zeta}
\def\O{{\cal O}}
\def\C{{\bf C}}
\def\R{{\bf R}}

\def\Re{{\rm Re}}
\def\Im{{\rm Im}}

\font\auth=cmcsc10
\def\jour{\sl}
\def\vol{\bf}

\topglue 1truein
\centerline{\bf A MORERA TYPE THEOREM IN THE STRIP}
\bigskip

\centerline{\auth A. Tumanov}
\bigskip

\beginsection  Introduction

We prove the following

{\bf Theorem.}
{\sl Let $f$ be a continuous function in the strip
$|\Im z|\le1$.
Suppose for every $r\in\R$ the restriction of $f$
to the circle $|z-r|=1$ extends holomorphically inside
the circle. Then $f$ is holomorphic in the strip
$|\Im z|<1$.}
\medskip

The result can be regarded as a Morera type theorem
because the holomorphic extendibility is equivalent to
the moment condition. For various versions of the Morera
theorem, see e. g. [ABC], [BZ].

The question answered by the above theorem has been
known for over a decade.
A general natural question is whether the
analyticity of a function can be tested by restricting
the function on (usually one parameter) families of
Jordan curves.

Agranovsky and Val'sky [AV] proved the result for any family
of curves invariant under Euclidean motions of the plane.

Globevnik [G] proved tests of analyticity in an annulus
for any rotation invariant family of Jordan curves.
Apparently, the compactness of the group of rotations
played an important role, because for the family obtained
by translating a given curve parallel to the real line,
the problem has since been open even for the circle.

Agranovsky [A] proved Theorem under additional growth
assumptions on $f$. He also proved the result in case
$f$ is real analytic in the disc $|z|<1+\epsilon$.

Agranovsky and Globevnik [AG] have recently solved
the problem for arbitrary one parameter families
of circles for
rational functions of two real variables and for
real analytic functions. In particular, they found
families on which the analyticity cannot be tested.

Ehrenpreis [E] has also proved Theorem for real
analytic $f$, but his results hold for PDE more general
than the Cauchy--Riemann equation.

Despite the one variable nature of the problem,
we use the analysis of several complex variables,
specifically, the extendibility of CR functions.
The proof goes through for a finite strip and
some curves other than the circle, but for a general
curve the question remains open.

The author thanks Mark Agranovsky for
useful discussions.

\beginsection  1. Proof of the main result

We prove the theorem stated in the introduction.

{\it Step 1.} We introduce a CR function $F$ on a
hypersurface $M\subset\C^2$.

Denote by $f_r(z)$ the holomorphic extension of $f$
to the disc $|z-r|<1$.
Define $F(z,w)=f_r(z)$ for $w=z-r$.
The function $F$ is defined on the real hyperplane
$\Im z=\Im w$ for $|w|\le 1$ and holomorphic on
every disc $z=w+r$, $|w|<1$, where $r\in\R$ is constant.
Hence $F$ is a CR function.
(For a reference on CR functions, see e. g. [BER].)

On the boundary $E=\{(z,w): \Im z=\Im w, |w|=1 \}$
of the region on the hyperplane where $F$ is defined,
we have $F(z,w)=f(z)$.
Note that $G(w)=-w^{-1}$ reflects the circle
$|w|=1$ about the imaginary axis.
Then $F(z,G(w))=f(z)$, too.
Hence $F$ satisfies $F(z,w)=F(z,G(w))$ on $E$.

Define the hypersurface $M=\{(z,w):\Im z=h(w), w\in\C \}$,
where $h(w)=\Im w$ for $|w|\le1$ and
$h(w)=h(G(w))=-\Im (w^{-1})$ for $|w|\ge1$.
Then $h$ is continuous and it is smooth
except for $|w|=1$.

Extend $F$ to the rest of $M$ by $F(z,w)=F(z,G(w))$.
Then $F$ is a continuous CR function on $M$.

We will ultimately prove that $F$ extends
holomorphically to
a neighborhood of $(0,0)$ and that actually
$F(z,w)$ is independent of $w$. This will mean
that $f_r(z)$ is independent of $r$, whence
$f$ is holomorphic.
\medskip

{\it Step 2.}
We prove that $F$ extends holomorphically to
a one-sided neighborhood of each point of $E$
except the points $(z,w)\in E$, $w=\pm 1$.

Indeed, $M$ consists of two smooth hypersurfaces
that meet transversally at the totally real
edge $E$ except the points where $w=\pm 1$
in which $E$ has complex tangencies.
By the ``edge of the wedge'' theorem by
Ayrapetian and Henkin [AH], in a neighborhood
of every point $(z_0,w_0)\in E$ such that
$\pm\Im w_0>0$, $F$ extends holomorphically
to the side $\pm(\Im z-h(w))<0$ respectively,
which completes this step.
\medskip

{\it Step 3.}
We prove that $F$ extends to a (full) neighborhood
of every point $(z,w)\in E$ with $w=\pm 1$.

This is done by extending $F$ into small analytic
discs. Let $\Delta$ be the standard unit disc in
complex plane. Recall an analytic disc in $\C^2$
is a map $g:\Delta\to\C^2$,
$g\in C(\bar\Delta)\cap\O(\Delta)$.
We say the disc $g$ is attached to $M$ if
$g(b\Delta)\subset M$.
By the Baouendi-Treves [BT] approximation theorem,
the CR function $F$ is locally a uniform limit of
holomorphic polynomials. Hence $F$ extends
holomorphically to every open set covered by
(the images of) small analytic discs attached
to $M$. (See e. g. [BER].)

Small analytic discs attached to $M$ of the form
$\Im z=h(w)$ are constructed explicitly as follows.
Let $\z\mapsto g(\z)=(z(\z)=x(\z)+iy(\z),w(\z))$,
$\z\in\Delta$, be an analytic disc attached to $M$.
Given the w-component $\z\mapsto w(\z)$ and
$x(0)=x_0$, we have $y(\z)=h(w(\z))$ for $|\z|=1$,
hence the z-component
$\z\mapsto z(\z)$ is defined as
the unique holomorphic function in $\Delta$ with given
imaginary part on the circle and given $x(0)$.
The ``center'' $g(0)=(x_0+iy_0,w(0))$ is found by
$$
y_0={1\over 2\pi i}\int_{|\z|=1}h(w(\z))\,{d\z\over\z}.
\eqno(1)
$$

We construct the family of discs
$\z\mapsto g(\z)=g(\z,x_0,w_0,t)$ for which
$w(\z)=w_0+a\z+t\phi(\z)$, where
$a>0$ is fixed and small enough to make
the discs small, and $\phi$ is a fixed holomorphic
function with $\phi(0)=0$, which we will choose later.
The discs depend on the parameters
$x_0\in\R$, $w_0\in\C$ close to $\pm1$, and $t\in\R$
close to 0.

Consider the evaluation map
$\Phi:(x_0,w_0,t)\mapsto g(0,x_0,w_0,t)$.
Note $\Phi(x_0,\pm1,0)=(x_0,\pm1)$.
We will prove that $\Phi$ is a diffeomorphism
in a neighborhood of $(x_0,\pm1,0)$, so
the discs cover a neighborhood of $(x_0,\pm1)$,
which will complete Step 3.

By the implicit function theorem, it suffices
to show that
$\dot y(0)={d\over dt}\big|_{t=0}y(0)\ne0$
for $w_0=\pm1$. Note that $h(w)=\Im H(w)$, where
$H(w)=w$ for $|w|\le1$ and
$H(w)=G(w)$ for $|w|>1$. Then (1) turns into
$$
y(0)=\Im \(
{1\over 2\pi i}
\int_{|\z|=1}H(\pm1+a\z+t\phi(\z))\,{d\z\over\z}
\).
\eqno(2)
$$
Differentiating (2) yields
$$
\dot y(0)=\Im \(
{1\over 2\pi i}
\int_{|\z|=1}H'(\pm1+a\z)\phi(\z)\,{d\z\over\z}
\).
$$
If the above expression vanished for all
polynomials $\phi$, then by the moment conditions,
$\z\mapsto H'(\pm1+a\z)$ would extend holomorphically
from $b\Delta$ to $\Delta$, which is not the case.
Hence $\dot y(0)\ne0$ for some $\phi$, which
completes Step 3.
\medskip

{\it Step 4.} We prove that $F$ extends holomorphically
to a neighborhood of $M\setminus E$.

Indeed, $M\setminus E$ is a union of complex curves,
each of which contains a boundary point of the form
$(x_0,1)$, where $x_0\in\R$.
By Step 3, $F$ is holomorphic at that point,
and the conclusion follows
by propagation of analyticity along complex curves
(Hanges and Treves [HT]).
\medskip

{\it Step 5.} We prove that $F(z,w)$ is independent of $w$.

We extend $F$ to ``big'' discs attached to $M$.
Extendibility to big discs does not follow
from the approximation theorem.
It follows by the continuity principle because
$F$ is holomorphic in a (one-sided near $E$)
neighborhood of $M$ and because discs attached
to $M$ can be contracted into a point.
We don't have to pay attention to the fact
that the neighborhood in which
$F$ is holomorphic is one-sided near $E$,
because $M$ can be slightly pushed inside
the neighborhood.

Consider the following one parameter family
of discs $D(a)$, $a>0$, with centers at $(x_0,0)$,
$x_0\in\R$:
$$
|w|<a,
\qquad
z=\left\{ \matrix{
x_0+w,\;\;\;\;\;\; & a\le1 \cr
x_0+a^{-2}w, & a>1 \cr}\right.
$$
We now regard the discs $D(a)$ as subsets in $\C^2$
rather than mappings of the standard disc.
The discs are attached to $M$;
they contract to $(x_0,0)$ as $a\to 0$.
The analytic continuations of $F$ along all
$D(a)$ match because the discs have the same
center.

We claim that $F|_{D(a)}$ is bounded uniformly
in $a$.
Indeed, by the construction in Step 1,
$F|_{bD(a)}$ for $a>1$ is defined using $f_r(z)$
with $r=z-G(w)$, where $|w|=a$, $z=x_0+a^{-2}w$.
Hence $r=x_0+a^{-2}w+w^{-1}=x_0+a^{-2}(w+\bar w)$,
$$
|r|\le|x_0|+2a^{-1},
\eqno(3)
$$
so $r$ whence $f_r$
is uniformly bounded in $a$.

Now for small $x_0$ and $\epsilon>0$,
$\phi(w)=F(x_0,w)$ is a uniform limit of
$\phi_a(w)=F(x_0+a^{-2}w,w)$ in $|w|<\epsilon$
as $a\to\infty$. The function $\phi_a$
is holomorphic in $|w|<a$ and is bounded
uniformly in $a$. Hence $\phi$ is constant.
This completes the last step in the proof of
the theorem.
\medskip

{\bf Remark.}
Theorem also holds in a finite strip.
Indeed, for Steps 1--4 it suffices to assume
the hypothesis of the theorem for
$|r|<1+\epsilon$ with arbitrarily small
$\epsilon>0$.
This seems to be a natural range for $r$ in
the theorem.
However, the contraction of the discs $D(a)$
in Step 5 as it follows from (3), requires
the range $|r|<2+\epsilon$.
One may try to reduce this range by
using a different contraction.

\beginsection  2. Proof for other curves

Our method applies to some curves other than
the circle.

{\bf Theorem'} {\sl
Let $a>0, b>0$.
Let $f$ be a continuous function in the strip
$|\Im z|\le b$.
Suppose for every $r\in\R$ the restriction of $f$
to the ellipse
$a^{-2}(x-r)^2+b^{-2}y^2=1$, where $z=x+iy$,
extends holomorphically inside
the ellipse. Then $f$ is holomorphic in the strip
$|\Im z|<b$.
}

{\it Proof.}
The proof consists of the same five steps.
We give a proof in the case $a<b$, the other case being
similar, even simpler. We normalize $a$ and $b$ so that
$b^2-a^2=1$, then we fix a number $p>0$, such that
$a=\sinh p$, $b=\cosh p$.
Denote by $E_q$ the ellipse with half-axes
$(\sinh q, \cosh q)$ centered at the origin and by
$D_q$ the domain bounded by $E_q$.

Following Step 1, we introduce $F(z,w)=f_r(z)$ for
$w=z-r$, where $f_r(z)$ is the extension of $f$
into $D_p+r$, $r\in\R$. Then $F$ is a CR function on
the hyperplane $\Im z=\Im w$, $w\in D_p$.

For $q\in\R$, we introduce the function
$$
G_q(w)=(\cosh q) w+(\sinh q) \sqrt{w^2+1}
$$
holomorphic in the exterior of the line segment
$[-i,i]$, where we choose the branch of the
square root such that $\sqrt{w^2+1}>0$ for $w>0$.
Note $G_q$ is an odd function.

For $q,r,t\in\R$,
$w=\sinh(q+it)=\sinh q \cos t+i\cosh q \sin t\in E_q$,
we have
$G_r(w)=\sinh(\pm r+q+it)$ if $\pm q>0$
respectively.
Hence $G_r(E_q)=E_{|r+q|}$ for $q>0$.

From the above properties it follows that
$G_{-2p}$ gives a conformal mapping
of $D_{2p}\setminus[-i,i]$ to itself, and
reflects $E_p$ about the imaginary axis.

We define $M$ by the equation $\Im z=h(w)$, where
$h(w)=\Im w$ for $w\in D_p$, and extends
to the whole plane applying the relation
$h(w)=h(G_{-2p}(w))$ infinitely many times.
We extend $F$ to the rest of $M$ so that
$F(z,w)=F(z,G_{-2p}(w))$ is also applied
infinitely many times.

This completes Step 1 of the proof.
Steps 2--4 go through along the same lines.

In Step 5 we consider the discs $D(n)$
given by $z=x_0+{w\over\cosh(2np)}$,
$w\in D_{2p}$, where $n$ is a positive integer.
Then for small $w$, the function
$\phi(w)=F(x_0,w)$ is approximated by
$\phi_n(w)=F(x_0+{w\over\cosh(2np)},w)$
as $n\to\infty$.
The function $\phi_n$ is defined in $D_{2np}$
and uniformly bounded in $n$.
Hence, $\phi$ is constant, and the
proof is complete.
\medskip

We consider one more example of a curve for
which our method works. This may be of interest
because the curve meets the boundaries of the
horizontal strip at a nonzero angle.
Let $\Gamma=\Gamma_+\cup\Gamma_-$, where
$\Gamma_\pm$ are the circular arcs defined as follows
$$
\Gamma_\pm=
\{z: |z\sin\alpha \pm \cos\alpha|=1, \pm\Re z\ge0 \}.
$$
The arcs intersect at the angle $2\alpha$,
where we choose $\alpha=\pi/n$, $n\ge 3$ is integer.
The proof is similar.
The defining function $h$ of the hypersurface $M$
is constructed using the relation $h(w)=h(G(w))$,
where
$G(w)=(w+\tan\alpha)(-w\tan\alpha+1)^{-1}$
maps $\Gamma_-$ to  $\Gamma_+$.

\beginsection References

\frenchspacing

\item{[A]} M. Agranovsky, Private communication.

\item{[AG]} M. Agranovsky and J. Globevnik,
Analyticity on circles for rational and real analytic
functions of two real variables,
Preprint series, Institute of Math., Phys., and Mech.,
University of Ljubljana,
vol. 39 (2001), 756, 1-43.

\item{[ABC]} M. Agranovsky, C. Berenstein,and D.-C. Chang,
Morera theorem for holomorphic $H^p$ spaces in the
Heisenberg group,
{\jour J. Reine und Angew. Math. \vol 443}
(1993), 49--89.

\item{[AV]} M. Agranovsky and R. Val'sky,
Maximality of invariant algebras of functions,
{\jour Siberian Math. J. \vol 12}
(1971), 1--7.

\item{[AH]}
R. A. Ayrapetian and G. M. Henkin,
Analytic continuation of CR functions through the
``edge of the wedge''.
{\jour Dokl. Akad. Nauk SSSR \vol 259}
(1981), 777--781.

\item{[BER]}
M. S. Baouendi, P. Ebenfelt, and L. P. Rothschild,
{\sl Real Submanifolds in Complex Space and their
Mappings}, Princeton Math. Series 47, Princeton
Univ. Press, 1999.

\item{[BT]}
M. S. Baouendi and F. Treves,
A property of the functions and distributions annihilated
by a locally integrable system of complex vector fields,
{\jour Ann. of Math. \vol 114}
(1981), 387--421.

\item{[BZ]}
C. Berenstein and L. Zalcman,
Pompeiu's problem on symmetric spaces,
{\jour Comm. Math. Helv. \vol 55}
(1980), 593-621.

\item{[E]}
L. Ehrenpreis,
Three problems at Mount Holyoke,
{\jour Contemp. Math., \vol 278},
Providence, RI, 2001.

\item{[G]}
J. Globevnik,
Testing analyticity on rotation invariant families of curves,
{\jour Trans. Amer. Math. Soc. \vol 306}
(1988), 401-410.

\item{[HT]}
N. Hanges and F. Treves,
Propagation of holomorphic extendibility of CR function,
{\jour Math. Ann. \vol 263}
(1983), 157--177.

\bigskip
\bigskip\noindent
Alexander Tumanov, Department of Mathematics, University of Illinois,
Urbana, IL 61801. E-mail: tumanov@math.uiuc.edu
\bye